# A DECOMPOSITION BY NON-NEGATIVE FUNCTIONS IN THE SOBOLEV SPACE $W^{k,1}$

AUGUSTO C. PONCE AND DANIEL SPECTOR

ABSTRACT. We show how a strong capacitary inequality can be used to give a decomposition of any function in the Sobolev space $W^{k,1}(\mathbb{R}^d)$ as the difference of two non-negative functions in the same space with control of their norms.

## 1. INTRODUCTION

In this paper we are interested in the following problem in the study of weakly differentiable functions:

*Question.* Given $u$ in the Sobolev space $W^{k,p}(\mathbb{R}^d)$ for $k \in \mathbb{N}_*$ and $p \geq 1$, is it possible to find non-negative functions $u^\oplus$ and $u^\ominus$ in the same space such that $u = u^\oplus - u^\ominus$ almost everywhere in $\mathbb{R}^d$ with control of the $W^{k,p}$ norms of $u^\oplus$ and $u^\ominus$?

An affirmative answer to this question has the practical consequence of enabling the qualification of arguments in $W^{k,p}(\mathbb{R}^d)$ with the statement "without loss of generality we assume that $u$ is non-negative, since $u$ can be decomposed into positive and negative parts...", which is useful, for example, in the proofs of various Hardy and Sobolev inequalities. This is standard in the first-order case for any $p \geq 1$, while with a little thought one can give a simple solution in the higher-order case for any $p > 1$. The main contribution of this work is the following theorem which gives such a decomposition when $p = 1$:

**Theorem 1.1.** *Let $k \in \mathbb{N}_*$. For every $u \in W^{k,1}(\mathbb{R}^d)$, there exist non-negative functions $u^\oplus, u^\ominus \in W^{k,1}(\mathbb{R}^d)$ such that*

$$u = u^\oplus - u^\ominus \quad \text{with} \quad \|u^\oplus\|_{W^{k,1}(\mathbb{R}^d)} + \|u^\ominus\|_{W^{k,1}(\mathbb{R}^d)} \leq C\|u\|_{W^{k,1}(\mathbb{R}^d)},$$

*for some constant $C > 0$ depending on $k$ and $d$.*

In the case of $W^{1,p}(\mathbb{R}^d)$, one can achieve such a decomposition for any $p \geq 1$, in the same manner as in $L^p(\mathbb{R}^d)$ and $C_0(\mathbb{R}^d)$, via the classical decomposition of $u$ as the difference between its positive and negative parts:

(1.1) $\qquad u^+ := \max\{u, 0\} \quad \text{and} \quad u^- := \max\{-u, 0\}.$

However such a choice is not suitable for higher order derivatives. For example, when $u \in W^{2,p}(\mathbb{R}^d)$ it may happen that the distributions $D^2 u^+$ and $D^2 u^-$ are merely finite measures with singular parts. This is not a simple matter of the lack of smoothness of the chosen truncation, as Dahlberg [3,5] has shown that if for a smooth function $H : \mathbb{R} \to \mathbb{R}$ one has

$$H(u) \in W^{2,p}(\mathbb{R}^d)$$

for all $u \in W^{2,p}(\mathbb{R}^d)$ with $1 < p < d/2$, then necessarily $H(t) = ct$ for some $c \in \mathbb{R}$. Such a conclusion propagates to $W^{k,p}(\mathbb{R}^d)$ for any $k \geq 3$ and $1 \leq p < d/k$.

As an alternative to composition, there is for $1 < p < \infty$ a convenient approach which relies on the characterization of $W^{k,p}(\mathbb{R}^d)$ via Bessel potentials: One has $u \in W^{k,p}(\mathbb{R}^d)$ if and only if $u = g_k * f$ for some $f \in L^p(\mathbb{R}^d)$ with

$$\|f\|_{L^p(\mathbb{R}^d)} \sim \|u\|_{W^{k,p}(\mathbb{R}^d)}.$$

Since $g_k \geq 0$, one can take the decomposition

$$u = g_k * f^+ - g_k * f^-,$$

for which one has the desired control of the norms.

There is no such an equivalence for $p = 1$, and in this regime one requires a replacement for harmonic analysis techniques with more geometric ideas. In particular, our proof of Theorem 1.1 relies on a strong capacitary inequality which has its roots on the Boxing inequality of Gustin [7]. The reader is perhaps more familiar with weak-type capacitary inequalities involving the Hardy-Littlewood maximal function, for example,

$$(1.2) \qquad \operatorname{cap}_{W^{1,1}}(\{\mathcal{M}u > t\}) \leq \frac{C}{t}\|u\|_{W^{1,1}(\mathbb{R}^d)},$$

for every $t > 0$ and $u \in W^{1,1}(\mathbb{R}^d)$. This inequality has been pioneered by Federer and Ziemer [6] and is a powerful tool for establishing fine properties of functions, see [15].

To express a strong-type analogue which is convenient for our purposes, we rely on the Choquet integral with respect to the $W^{k,1}$ capacity, defined in analogy with Cavalieri's principle for measures as

$$\int_{\mathbb{R}^d} |\varphi| \, \mathrm{d}\operatorname{cap}_{W^{k,1}} := \int_0^\infty \operatorname{cap}_{W^{k,1}}(\{|\varphi| > t\}) \, \mathrm{d}t.$$

Then the strong capacitary inequality we use to prove Theorem 1.1 is our

**Theorem 1.2.** *Let $k \in \mathbb{N}_*$. For every $\varphi \in C_c^\infty(\mathbb{R}^d)$, we have*

$$\int_{\mathbb{R}^d} |\varphi| \, \mathrm{d}\operatorname{cap}_{W^{k,1}} \leq C'\|\varphi\|_{W^{k,1}(\mathbb{R}^d)},$$

*for some constant $C' > 0$ depending on $k$ and $d$.*

Such inequalities have been initiated in the first-order case by Maz'ya [9] and later pursued for any order and $p > 1$ in [1, 10] (see also Theorem 1.1.2 in [11] for $k \leq d$). Our result completes the picture in the regime $p = 1$, as in fact we prove a more general result than Theorem 1.2 which includes the scale of fractional Sobolev spaces $W^{\alpha,1}$ for any $\alpha > 0$, defined in terms of the Gagliardo semi-norm. Since Theorem 1.1 is directly argued from such an inequality, we obtain a decomposition that is also valid on the $W^{\alpha,1}$ scale, see Theorem 3.1 below. In contrast to (1.2), one cannot replace $\varphi$ in the left-hand side with the maximal function $\mathcal{M}\varphi$. Nonetheless, as we show in Section 4, one does have a strong-type inequality with a local variant of the maximal function, which contains both Theorem 1.2 and the inequality (1.2). Let us finally remark that Theorem 1.2 is most interesting in the range $0 < \alpha \leq d$, as we explain below that all the Sobolev capacities are equivalent for $\alpha \geq d$.

The plan of the paper is as follows. In Section 2 we show how recent work of the authors [13] on the Boxing inequality of Gustin implies a non-homogeneous form of this inequality that involves the Choquet integral with respect to the Hausdorff outer measures $\mathcal{H}_\delta^{d-\alpha}$ for $0 < \alpha \leq d$ and $0 < \delta < \infty$. In this range of $\alpha$ such an estimate is equivalent to the strong capacitary inequality in Theorem 1.2. We then argue the case $\alpha \geq d$ from the observation that all capacities are equivalent to $\mathcal{H}_\delta^0$, see Proposition 2.4. In Section 3 we show how the strong capacitary inequality implies Theorem 1.1 and its fractional counterpart. In Section 4, we prove the local maximal-function counterpart of the strong capacitary inequality for functions in $W^{\alpha,1}(\mathbb{R}^d)$ that implies both (1.2) and Theorem 1.2. In Section 5 we rely on a lemma of Harvey and Polking's [8] to show that when $u$ is bounded, the functions $u^\oplus$ and $u^\ominus$ can inherit the same property.

## 2. Non-homogeneous Boxing inequality

Let $k \in \mathbb{N}$ and denote by $W^{k,1}(\mathbb{R}^d)$ the space of weakly differentiable functions $u \in L^1(\mathbb{R}^d)$ such that $D^i u \in L^1(\mathbb{R}^d, \mathbb{R}^{d^i})$ for $i \in \{1, \ldots, k\}$. We define the norm of $u$ by

$$\|u\|_{W^{k,1}(\mathbb{R}^d)} := \sum_{i=0}^{k} \|D^i u\|_{L^1(\mathbb{R}^d)},$$

where $D^0 u := u$. Next, for $\alpha > 0$ non-integer, write $\alpha = k + \theta$ where $k \in \mathbb{N}$ and $0 < \theta < 1$. We say that $u \in W^{\alpha,1}(\mathbb{R}^d)$ whenever $u \in W^{k,1}(\mathbb{R}^d)$ and

$$\|u\|_{W^{\alpha,1}(\mathbb{R}^d)} := \|u\|_{W^{k,1}(\mathbb{R}^d)} + [D^k u]_{W^{\theta,1}(\mathbb{R}^d)}$$

is finite, where

$$[D^k u]_{W^{\theta,1}(\mathbb{R}^d)} := \int_{\mathbb{R}^d} \int_{\mathbb{R}^d} \frac{|D^k u(x) - D^k u(y)|}{|x-y|^{\theta+d}} \,\mathrm{d}x\,\mathrm{d}y$$

is the Gagliardo semi-norm of order $\theta$ of $D^k u$.

We rely on the following non-homogeneous version of Gustin's Boxing inequality:

**Proposition 2.1.** *Let $0 < \alpha \le d$ and $0 < \delta < \infty$. For every $\varphi \in C_c^\infty(\mathbb{R}^d)$, we have*

$$\int_{\mathbb{R}^d} |\varphi| \, \mathrm{d}\mathcal{H}_\delta^{d-\alpha} \le C \|\varphi\|_{W^{\alpha,1}(\mathbb{R}^d)},$$

*for some constant $C > 0$ depending on $\alpha$, $d$ and $\delta$.*

For every $0 < \delta \le \infty$, $\mathcal{H}_\delta^{d-\alpha}$ is the Hausdorff outer measure defined for every $A \subset \mathbb{R}^d$ by

$$\mathcal{H}_\delta^{d-\alpha}(A) := \inf \left\{ \sum_{i=1}^\infty r_i^{d-\alpha} : A \subset \bigcup_{i=1}^\infty B_{r_i}(x_i), \ r_i \le \delta \right\}$$

and the Choquet integral of $|\varphi|$ with respect to $\mathcal{H}_\delta^{d-\alpha}$ is

$$\int_{\mathbb{R}^d} |\varphi| \, \mathrm{d}\mathcal{H}_\delta^{d-\alpha} := \int_0^\infty \mathcal{H}_\delta^{d-\alpha}(\{|\varphi| > t\}) \, \mathrm{d}t.$$

In this definition one can replace the sets $\{|\varphi| > t\}$ by $\{|\varphi| \ge t\}$ without changing the value of the integral.

We prove Proposition 2.1 using the following strong form of the Boxing inequality for the Hausforff content $\mathcal{H}_\infty^{d-\alpha}$:

(2.1) $$\int_0^\infty |\varphi| \, \mathrm{d}\mathcal{H}_\infty^{d-\alpha} \le C'[\varphi]_{W^{\alpha,1}(\mathbb{R}^d)}, \quad \text{for every } \varphi \in C_c^\infty(\mathbb{R}^d).$$

The homogeneous semi-norm in the right-hand side is

$$[\varphi]_{W^{\alpha,1}(\mathbb{R}^d)} := \|D^k \varphi\|_{L^1(\mathbb{R}^d)},$$

when $\alpha = k \in \mathbb{N}$ and

$$[\varphi]_{W^{\alpha,1}(\mathbb{R}^d)} := [D^k \varphi]_{W^{\theta,1}(\mathbb{R}^d)},$$

when $\alpha \notin \mathbb{N}$ and $\alpha = k + \theta$ with $k \in \mathbb{N}$ and $0 < \theta < 1$. We refer the reader to [2, 13] for the proof of (2.1) when $\alpha < d$. The case $\alpha = d$ is straightforward: From the facts that $\mathcal{H}_\infty^0(\emptyset) = 0$ and $\mathcal{H}_\infty^0(A) = 1$ for every nonempty bounded subset $A \subset \mathbb{R}^d$ one has

$$\int_0^\infty |\varphi| \, \mathrm{d}\mathcal{H}_\infty^0 = \int_0^{\|\varphi\|_{L^\infty(\mathbb{R}^d)}} \mathcal{H}_\infty^0(\{|\varphi| > t\}) \, \mathrm{d}t = \|\varphi\|_{L^\infty(\mathbb{R}^d)}$$

and then (2.1) with $\alpha = d$ is equivalent to the classical inequality

(2.2) $$\|\varphi\|_{L^\infty(\mathbb{R}^d)} \le C' \|D^d \varphi\|_{L^1(\mathbb{R}^d)}.$$

We deduce Proposition 2.1 for $\alpha < d$ using the next

**Lemma 2.2.** *Let $(Q_j)_{j\in\mathbb{N}}$ be a family of closed cubes in $\mathbb{R}^d$ with disjoint interiors obtained by translation of a fixed cube as $Q_j = Q + c_j$, and let $\zeta_j : \mathbb{R}^d \to \mathbb{R}$ be defined by $\zeta_j(x) = \zeta(x - c_j)$ for some $\zeta \in C_c^\infty(2Q)$. Given $\alpha > 0$, there exists a constant $C'' > 0$ depending on $\alpha$, $d$ and the side length of Q such that*

$$\sum_{j=0}^{\infty} \|u\zeta_j\|_{W^{\alpha,1}(\mathbb{R}^d)} \leq C'' \|u\|_{W^{\alpha,1}(\mathbb{R}^d)}, \quad \text{for every } u \in W^{\alpha,1}(\mathbb{R}^d).$$

*Proof of Lemma 2.2.* The estimate of the terms $\sum_{j=0}^{\infty} \|D^i(u\zeta_j)\|_{L^1(\mathbb{R}^d)}$ with $i \in \{0,\ldots,k\}$ is straightforward due to the local character of the $L^1$ norm. We thus assume that $\alpha = k+\theta$ with $0 < \theta < 1$ and estimate $\sum_{j=0}^{\infty} [D^k(u\zeta_j)]_{W^{\theta,1}(\mathbb{R}^d)}$. By the explicit formula of $D^k(u\zeta_j)$, it suffices to estimate the Gagliardo seminorm of the functions $D^i u \otimes D^{k-i}\zeta_j$. To simplify the notation we perform the estimates using $f := D^i u$ and $g_j := D^{k-i}\zeta_j$. We then have

$$[D^i u \otimes D^{k-i}\zeta_j]_{W^{\theta,1}(\mathbb{R}^d)} \leq \int_{\mathbb{R}^d} |g_j(x)| \left( \int_{\mathbb{R}^d} \frac{|f(x) - f(y)|}{|x-y|^{\theta+d}} \, dy \right) dx$$
$$+ \int_{\mathbb{R}^d} |f(y)| \left( \int_{\mathbb{R}^d} \frac{|g_j(x) - g_j(y)|}{|x-y|^{\theta+d}} \, dx \right) dy.$$

Since $g_j$ is supported in $2Q_j$ and the number of overlaps of these cubes is uniformly bounded by some constant depending on the dimension $d$, we have

$$\sum_{j=0}^{\infty} |g_j(x)| \leq C_1, \quad \text{for every } x \in \mathbb{R}^d.$$

For the second term, take $y \in \mathbb{R}^d$ and let $J_y \subset \mathbb{N}$ be the set of indices $j$ such that $y \in 4Q_j$. The number of elements in $J_y$ is bounded from above independently of $y$. By the uniform boundedness of $(g_j)_{j\in\mathbb{N}}$ and $(Dg_j)_{j\in\mathbb{N}}$, we then get

$$\sum_{j\in J_y} \int_{\mathbb{R}^d} \frac{|g_j(x) - g_j(y)|}{|x-y|^{\theta+d}} \, dx \leq C_2, \quad \text{for every } y \in \mathbb{R}^d.$$

Denoting by $\eta > 0$ the side length of $Q$, for $j \in \mathbb{N} \setminus J_y$ we have $d(y, 2Q_j) \geq \eta$ and also $g_j(y) = 0$. Thus,

$$\sum_{j\in\mathbb{N}\setminus J_y} \int_{\mathbb{R}^d} \frac{|g_j(x) - g_j(y)|}{|x-y|^{\theta+d}} \, dx \leq \|g_j\|_{L^\infty(\mathbb{R}^d)} \sum_{j\in\mathbb{N}\setminus J_y} \int_{2Q_j} \frac{dx}{|x-y|^{\theta+d}}$$
$$\leq C_3 \int_{\mathbb{R}^d \setminus B_\eta(y)} \frac{dx}{|x-y|^{\theta+d}} \leq C_4,$$

for every $y \in \mathbb{R}^d$. Hence,

$$\sum_{j=0}^{\infty} [D^i u \otimes D^{k-i}\zeta_j]_{W^{\theta,1}(\mathbb{R}^d)} \leq C_1 [f]_{W^{\theta,1}(\mathbb{R}^d)} + (C_2+C_4)\|f\|_{L^1(\mathbb{R}^d)} \leq C_5 \|u\|_{W^{\alpha,1}(\mathbb{R}^d)},$$

which implies the estimate of $\sum_{j=0}^{\infty} [D^k(u\zeta_j)]_{W^{\theta,1}(\mathbb{R}^d)}$. □

*Proof of Proposition 2.1.* Let $Q$ be the cube centered at the origin with diameter $\delta/2$. We cover $\mathbb{R}^d$ using a countable collection of closed cubes $(Q_j)_{j \in \mathbb{N}}$ with disjoint interiors of the form $Q_j = Q + c_j$. By subadditivity of $\mathcal{H}_\delta^{d-\alpha}$, for every subset $A \subset \mathbb{R}^d$ we have

$$(2.3) \quad \mathcal{H}_\delta^{d-\alpha}(A) \leq \sum_{j=0}^{\infty} \mathcal{H}_\delta^{d-\alpha}(A \cap Q_j) = \sum_{j=0}^{\infty} \mathcal{H}_\infty^{d-\alpha}(A \cap Q_j),$$

where the equality follows from the fact that each set $A \cap Q_j$ has diameter less than $\delta$. Let $\zeta \in C_c^\infty(2Q)$ be such that $\zeta = 1$ on $Q$, and consider the function $\zeta_j : \mathbb{R}^d \to \mathbb{R}$ defined by $\zeta_j(x) = \zeta(x - c_j)$. Observe that for every $\varphi \in C_c^\infty(\mathbb{R}^d)$ and $j \in \mathbb{N}$ we have

$$\{|\varphi| \geq t\} \cap Q_j \subset \{|\varphi\zeta_j| \geq t\}.$$

Applying the homogeneous Boxing inequality (2.1) to $\varphi\zeta_j$, we have

$$(2.4) \quad \int_0^\infty \mathcal{H}_\infty^{d-\alpha}(\{|\varphi\zeta_j| \geq t\}) \, dt \leq C'[\varphi\zeta_j]_{W^{\alpha,1}(\mathbb{R}^d)}.$$

A combination of (2.3) with $A = \{|\varphi| \geq t\}$ and (2.4) then gives

$$\int_0^\infty \mathcal{H}_\delta^{d-\alpha}(\{|\varphi| \geq t\}) \, dt \leq \sum_{j=0}^{\infty} \int_0^\infty \mathcal{H}_\infty^{d-\alpha}(\{|\varphi\zeta_j| \geq t\}) \, dt \leq C' \sum_{j=0}^{\infty} [\varphi\zeta_j]_{W^{\alpha,1}(\mathbb{R}^d)}.$$

By Lemma 2.2 we get

$$\int_0^\infty \mathcal{H}_\delta^{d-\alpha}(\{|\varphi| \geq t\}) \, dt \leq C_6 \|\varphi\|_{W^{\alpha,1}(\mathbb{R}^d)}. \quad \Box$$

One deduces from Proposition 2.1 and a straightforward adaptation of the proof of Theorem 2.1 in [13] the following:

**Corollary 2.3.** *For every $0 < \alpha \leq d$ and $0 < \delta < \infty$, we have*

$$\mathrm{cap}_{W^{\alpha,1}} \sim \mathcal{H}_\delta^{d-\alpha}.$$

Given $\alpha > 0$, the Sobolev capacity $\mathrm{cap}_{W^{\alpha,1}}$ is defined for every compact set $K \subset \mathbb{R}^d$ by

$$\mathrm{cap}_{W^{\alpha,1}}(K) := \inf \left\{ \|\varphi\|_{W^{\alpha,1}(\mathbb{R}^d)} : \varphi \in C_c^\infty(\mathbb{R}^d),\ \varphi \geq 0 \text{ in } \mathbb{R}^d,\ \varphi > 1 \text{ on } K \right\}.$$

It is then extended to the class of open sets $\omega \subset \mathbb{R}^d$ as the supremum

$$\mathrm{cap}_{W^{\alpha,1}}(\omega) := \sup \left\{ \mathrm{cap}_{W^{\alpha,1}}(K) : K \subset \omega \text{ compact} \right\}.$$

The equivalence in Corollary 2.3 has been proved by Carlsson and Maz'ya for $\alpha \in \mathbb{N}$ and $\alpha < d$, see Lemma 3 in [4]. On balls $B_r(x) \subset \mathbb{R}^d$, such a result can be easily obtained from a scaling argument. Indeed, when $\alpha \leq d$ a separate analysis in the regimes $r \leq 1$ and $r \geq 1$ yields

$$(2.5) \quad \mathrm{cap}_{W^{\alpha,1}}(B_r(x)) \sim \max\{r^{d-\alpha}, r^d\},$$

which is the same behavior as for $\mathcal{H}_\delta^{d-\alpha}(B_r(x))$. In the range $\alpha \geq d$, one uses the classical inequality (2.2) to get

(2.6) $$\operatorname{cap}_{W^{\alpha,1}}(B_r(x)) \sim \max\{1, r^d\}.$$

While the $W^{\alpha,1}$ capacity strongly depends on $\alpha$ for $\alpha \leq d$, we here observe that all $W^{\alpha,1}$ capacities are equivalent for $\alpha \geq d$:

**Proposition 2.4.** *For every $\alpha \geq d$ and $0 < \delta < \infty$, we have*
$$\operatorname{cap}_{W^{\alpha,1}} \sim \mathcal{H}_\delta^0.$$

*Proof.* Given a compact subset $K \subset \mathbb{R}^d$ and a finite covering $K \subset \bigcup_{i=1}^N B_{r_i}(x_i)$ with $r_i \leq \delta$, by subadditivbity of the capacity and (2.6),
$$\operatorname{cap}_{W^{\alpha,1}}(K) \leq \sum_{i=1}^N \operatorname{cap}_{W^{\alpha,1}}(B_{r_i}(x_i)) \leq \sum_{i=1}^N C_1 \max\{1, r_i^d\} \leq C_1 \max\{1, \delta^d\} N.$$
Minimizing the right-hand side with respect to the number of balls $N$, we get
$$\operatorname{cap}_{W^{\alpha,1}}(K) \leq C_1 \max\{1, \delta^d\} \mathcal{H}_\delta^0(K).$$
Applying Corollary 2.3 with $\alpha = d$, we get the reverse inequality
$$\mathcal{H}_\delta^0(K) \leq C_2 \operatorname{cap}_{W^{d,1}}(K) \leq C_2 \operatorname{cap}_{W^{\alpha,1}}(K). \qquad \square$$

Using Propositions 2.1 and 2.4, we now prove the strong capacitary inequality for the $W^{\alpha,1}$ norm, which includes Theorem 1.2 when $\alpha$ is integer:

**Theorem 2.5.** *Let $\alpha > 0$. For every $\varphi \in C_c^\infty(\mathbb{R}^d)$, we have*
$$\int_{\mathbb{R}^d} |\varphi| \,\mathrm{d}\operatorname{cap}_{W^{\alpha,1}} \leq C' \|\varphi\|_{W^{\alpha,1}(\mathbb{R}^d)},$$
*for some constant $C' > 0$ depending on $\alpha$ and $d$.*

*Proof.* The case of order $\alpha \leq d$ is contained in Proposition 2.1 and the straightforward inequality $\operatorname{cap}_{W^{\alpha,1}} \leq C' \mathcal{H}_\delta^{d-\alpha}$. When $\alpha > d$, one applies Proposition 2.4 and Proposition 2.1 with order $d$ to get
$$\int_0^\infty \operatorname{cap}_{W^{\alpha,1}}(\{|\varphi| > t\}) \,\mathrm{d}t \leq C_1 \int_0^\infty \mathcal{H}_1^0(\{|\varphi| > t\}) \,\mathrm{d}t$$
$$\leq C_2 \|\varphi\|_{W^{d,1}(\mathbb{R}^d)} \leq C_2 \|\varphi\|_{W^{\alpha,1}(\mathbb{R}^d)}. \qquad \square$$

## 3. Proof of Theorem 1.1

We establish the following theorem that includes Theorem 1.1 for integer orders:

**Theorem 3.1.** *Let $\alpha > 0$. For every $u \in W^{\alpha,1}(\mathbb{R}^d)$, there exist non-negative functions $u^\oplus, u^\ominus \in W^{\alpha,1}(\mathbb{R}^d)$ such that*
$$u = u^\oplus - u^\ominus \quad \text{with} \quad \|u^\oplus\|_{W^{\alpha,1}(\mathbb{R}^d)} + \|u^\ominus\|_{W^{\alpha,1}(\mathbb{R}^d)} \leq C \|u\|_{W^{\alpha,1}(\mathbb{R}^d)},$$
*for some constant $C > 0$ depending on $\alpha$ and $d$.*

Let us begin with an analogous property for smooth functions:

**Proposition 3.2.** *Let $\alpha > 0$. For every $\varphi \in C_c^\infty(\mathbb{R}^d)$, there exists $\psi \in C_c^\infty(\mathbb{R}^d)$ such that $|\varphi| \le \psi$ in $\mathbb{R}^d$ with*

$$\|\psi\|_{W^{\alpha,1}(\mathbb{R}^d)} \le C'\|\varphi\|_{W^{\alpha,1}(\mathbb{R}^d)},$$

*for some constant $C' > 0$ depending on $\alpha$ and $d$.*

For $\alpha \le 2$ this proposition can be proved by regularization of $|\varphi|$. In contrast, our argument works for every order $\alpha > 0$ and ultimately relies on the Hahn-Banach theorem through the strong form of the Boxing inequality (see e.g. the proof of Lemma 4.6 in [13]).

*Proof of Proposition 3.2.* We assume that $\varphi \not\equiv 0$. By the monotonicity of the capacity,

$$(3.1) \qquad \sum_{j=-\infty}^{\infty} 2^{j-1} \operatorname{cap}_{W^{\alpha,1}} (\{|\varphi| \ge 2^j\}) \le \int_0^\infty \operatorname{cap}_{W^{\alpha,1}} (\{|\varphi| \ge t\}) \, dt.$$

Let $J \in \mathbb{Z}$ be an integer such that

$$|\varphi| \le 2^J \quad \text{in } \mathbb{R}^d.$$

For every integer $j \le J$, take a non-negative function $\psi_j \in C_c^\infty(\mathbb{R}^d)$ such that $\psi_j > 1$ in $\{|\varphi| \ge 2^j\}$ and

$$\|\psi_j\|_{W^{\alpha,1}(\mathbb{R}^d)} \le 2 \operatorname{cap}_{W^{\alpha,1}} (\{|\varphi| \ge 2^j\}).$$

Let $m$ be an integer to be explicitly chosen below, depending on $\varphi$, and let $\zeta \in C_c^\infty(\mathbb{R}^d)$ be such that $0 \le \zeta \le 1$ in $\mathbb{R}^d$ and $\zeta = 1$ in $\operatorname{supp} \varphi$.

We claim that the function

$$\psi := 2^m \zeta + \sum_{j=m}^{J-1} 2^{j+1} \psi_j$$

satisfies

$$(3.2) \qquad \psi \ge |\varphi| \quad \text{in } \mathbb{R}^d.$$

We first observe that $\psi \ge 2^m$ in $\operatorname{supp} \varphi$. We thus have to prove (3.2) on the set $\{|\varphi| \ge 2^m\}$. To this end, given $x \in \mathbb{R}^d$ such that $2^m \le |\varphi(x)| < 2^J$, take $i \in \mathbb{Z}$ with $m \le i \le J-1$ such that $2^i \le |\varphi(x)| < 2^{i+1}$. Since all functions $\psi_j$ are non-negative and $\psi_i(x) > 1$, we have

$$\sum_{j=m}^{J-1} 2^{j+1} \psi_j(x) \ge 2^{i+1} \psi_i(x) > 2^{i+1} > |\varphi(x)|.$$

That is,

$$|\varphi| \le \sum_{j=m}^{J-1} 2^{j+1} \psi_j \quad \text{in } \{2^m \le |\varphi| < 2^J\},$$

and so also in the larger set $\{2^m \leq |\varphi| \leq 2^J\} = \{|\varphi| \geq 2^m\}$ by continuity of $\varphi$. Hence, (3.2) holds in $\operatorname{supp}\varphi$ and then, by nonnegativity of $\psi$, this pointwise inequality holds in the entire space $\mathbb{R}^d$.

We now claim that

$$(3.3) \quad \|\psi\|_{W^{\alpha,1}(\mathbb{R}^d)} \leq 2^m \|\zeta\|_{W^{\alpha,1}(\mathbb{R}^d)} + 8C' \|\varphi\|_{W^{\alpha,1}(\mathbb{R}^d)},$$

where $C' > 0$ is the constant in Theorem 2.5. Indeed, by the choice of $\psi_j$ and (3.1),

$$\|\psi\|_{W^{\alpha,1}(\mathbb{R}^d)} \leq 2^m \|\zeta\|_{W^{\alpha,1}(\mathbb{R}^d)} + \sum_{j=m}^{J-1} 2^{j+1} \|\psi_j\|_{W^{\alpha,1}(\mathbb{R}^d)}$$

$$\leq 2^m \|\zeta\|_{W^{\alpha,1}(\mathbb{R}^d)} + \sum_{j=m}^{J-1} 2^{j+2} \operatorname{cap}_{W^{\alpha,1}}(\{|\varphi| \geq 2^j\})$$

$$\leq 2^m \|\zeta\|_{W^{\alpha,1}(\mathbb{R}^d)} + 8 \int_0^\infty \operatorname{cap}_{W^{\alpha,1}}(\{|\varphi| \geq t\}) \, \mathrm{d}t.$$

Estimate (3.3) thus follows from the strong capacitary inequality.

To conclude the proof, it now suffices to choose $m \in \mathbb{Z}$ such that $2^m \|\zeta\|_{W^{\alpha,1}(\mathbb{R}^d)} \leq \|\varphi\|_{W^{\alpha,1}(\mathbb{R}^d)}$, which is possible since $\zeta$ is independent of $m$. $\square$

*Proof of Theorem 3.1.* Given $u \in W^{\alpha,1}(\mathbb{R}^d)$, by density of $C_c^\infty(\mathbb{R}^d)$ in this space we can write $u$ as a strongly convergent series $u = \sum_{j=0}^\infty \varphi_j$, where $\varphi_j \in C_c^\infty(\mathbb{R}^d)$ and

$$\sum_{j=0}^\infty \|\varphi_j\|_{W^{\alpha,1}(\mathbb{R}^d)} \leq 2\|u\|_{W^{\alpha,1}(\mathbb{R}^d)}.$$

Take $\psi_j \in C_c^\infty(\mathbb{R}^d)$ such that $|\varphi_j| \leq \psi_j$ in $\mathbb{R}^d$ and $\|\psi_j\|_{W^{\alpha,1}(\mathbb{R}^d)} \leq C_1 \|\varphi_j\|_{W^{\alpha,1}(\mathbb{R}^d)}$ as in Proposition 3.2. The function $u^\oplus := \sum_{j=0}^\infty \psi_j$ is non-negative in $\mathbb{R}^d$ and satisfies

$$\|u^\oplus\|_{W^{\alpha,1}(\mathbb{R}^d)} \leq C_1 \sum_{j=0}^\infty \|\varphi_j\|_{W^{\alpha,1}(\mathbb{R}^d)} \leq 2C_1 \|u\|_{W^{\alpha,1}(\mathbb{R}^d)}.$$

In particular, $u^\oplus$ belongs to $W^{\alpha,1}(\mathbb{R}^d)$. The conclusion readily follows with $u^\ominus := u^\oplus - u$, which is non-negative by the choice of $u^\oplus$. $\square$

**Remark 3.3.** The proof of Theorem 3.1 shows that for every $u \in W^{\alpha,1}(\mathbb{R}^d)$ one has

$$|u| \leq U \quad \text{almost everywhere in } \mathbb{R}^d,$$

where $U := \sum_{j=0}^\infty \psi_j$ is the nondecreasing limit of a sequence of non-negative functions in $C_c^\infty(\mathbb{R}^d)$ and satisfies the estimate

$$\|U\|_{W^{\alpha,1}(\mathbb{R}^d)} \leq C \|u\|_{W^{\alpha,1}(\mathbb{R}^d)},$$

with a constant $C > 0$ depending on $\alpha$ and $d$. Such a construction can be interpreted as a regularized replacement of the absolute-value function that is adapted to the space $W^{\alpha,1}(\mathbb{R}^d)$.

## 4. MAXIMAL STRONG CAPACITARY INEQUALITY

We now give an improvement of the strong capacitary inequality that involves the maximal operator and implies the weak capacitary inequality (1.2):

**Proposition 4.1.** *Let $\alpha > 0$ and $0 < \delta < \infty$. For every $u \in W^{\alpha,1}(\mathbb{R}^d)$, we have*

$$\int_{\mathbb{R}^d} \mathcal{M}_\delta u \, \mathrm{d}\operatorname{cap}_{W^{\alpha,1}} \leq C\|u\|_{W^{\alpha,1}(\mathbb{R}^d)},$$

*for some constant $C > 0$ depending on $\alpha$, $d$ and $\delta$.*

Here, $\mathcal{M}_\delta u : \mathbb{R}^d \to [0, +\infty]$ denotes the local maximal function

$$\mathcal{M}_\delta u(x) := \sup\left\{ \fint_{B_r(x)} |u| : 0 < r \leq \delta \right\}.$$

Proposition 4.1 relies on the inequality

(4.1) $$\int_{\mathbb{R}^d} \mathcal{M} u \, \mathrm{d}\mathcal{H}^{d-\alpha}_\infty \leq C'[u]_{W^{\alpha,1}(\mathbb{R}^d)},$$

for every $u \in W^{\alpha,1}(\mathbb{R}^d)$ and $0 < \alpha \leq d$, which is proved in [13] using the Boxing inequality (2.1) and D. Adams' maximal estimate for the Choquet integral [2, 12]

$$\int_{\mathbb{R}^d} \mathcal{M} u \, \mathrm{d}\mathcal{H}^{d-\alpha}_\infty \leq C'' \int_{\mathbb{R}^d} |u| \, \mathrm{d}\mathcal{H}^{d-\alpha}_\infty.$$

Observe that for $\alpha = d$ this inequality is equivalent to the straightforward

$$\|\mathcal{M} u\|_{L^\infty(\mathbb{R}^d)} \leq \|u\|_{L^\infty(\mathbb{R}^d)}.$$

*Proof of Proposition 4.1.* We first assume that $\alpha \leq d$. As in the statement of Lemma 2.2, we cover $\mathbb{R}^d$ with cubes $Q_j$ centered at $c_j$ that are obtained by translation from a fixed cube $Q$ and we assume that each $Q_j$ has side length $4\delta$. Let $\zeta \in C_c^\infty(2Q)$ be such that $\zeta = 1$ on $\frac{3}{2}Q$, and consider the function $\zeta_j : \mathbb{R}^d \to \mathbb{R}$ defined by $\zeta_j(x) = \zeta(x - c_j)$. Observe that if $x \in Q_j$, then $B_\delta(x) \subset \frac{3}{2}Q_j$. Thus, for every $t > 0$,

$$\{\mathcal{M}_\delta u > t\} \cap Q_j \subset \{\mathcal{M}(u\zeta_j) > t\}.$$

Applying (2.3) with $A = \{\mathcal{M}_\delta u > t\}$ and using the subadditivity and monotonicity of $\mathcal{H}^{d-\alpha}_\infty$, we get

$$\mathcal{H}^{d-\alpha}_\delta(\{\mathcal{M}_\delta u > t\}) \leq \sum_{j=0}^\infty \mathcal{H}^{d-\alpha}_\infty(\{\mathcal{M}_\delta u > t\} \cap Q_j) \leq \sum_{j=0}^\infty \mathcal{H}^{d-\alpha}_\infty(\{\mathcal{M}(u\zeta_j) > t\}).$$

We now integrate this estimate with respect to $t$ over $(0, \infty)$. By (4.1) applied to each function $u\zeta_j$ we thus have

$$\int_{\mathbb{R}^d} \mathcal{M}_\delta u \, d\mathcal{H}^{d-\alpha}_\delta \leq C' \sum_{j=0}^{\infty} \|u\zeta_j\|_{W^{\alpha,1}(\mathbb{R}^d)}.$$

The conclusion follows for $\alpha \leq d$ since $\mathcal{H}^{d-\alpha}_\delta \sim \text{cap}_{W^{\alpha,1}}$ (by Corollary 2.3) and the series in the right-hand side is bounded from above by $\|u\|_{W^{\alpha,1}(\mathbb{R}^d)}$ (by Lemma 2.2).

When $\alpha > d$, it suffices to apply Proposition 2.4 and the inequality at order $d$ to get

$$\int_{\mathbb{R}^d} \mathcal{M}_\delta u \, d\,\text{cap}_{W^{\alpha,1}} \leq C_2 \int_{\mathbb{R}^d} \mathcal{M}_\delta u \, d\mathcal{H}^0_\delta \leq C_3 \|u\|_{W^{d,1}(\mathbb{R}^d)} \leq C_3 \|u\|_{W^{\alpha,1}(\mathbb{R}^d)}.$$
□

The counterpart of Proposition 4.1 for the usual Hardy-Littlewood maximal function $\mathcal{M}u = \mathcal{M}_\infty u$ is false due to the same obstruction as for the strong $L^1$ maximal inequality: Any $\varphi \in C^\infty_c(\mathbb{R}^d)$ such that $\varphi(0) \neq 0$ satisfies

(4.2) $$\int_{\mathbb{R}^d} \mathcal{M}\varphi \, d\,\text{cap}_{W^{\alpha,1}} = +\infty.$$

To this end, we may assume that $\varphi(0) = 1$. Since $\mathcal{M}\varphi(x)$ is bounded from below by $C_1/|x|^d$ for large values of $|x|$, for $0 < t < 1/2$ one has

$$\{\mathcal{M}\varphi > t\} \supset \overline{B_{\epsilon/t^{1/d}}(0)},$$

for some $\epsilon > 0$. On the other hand, for every $r > 0$ and $x \in \mathbb{R}^d$,

(4.3) $$\text{cap}_{W^{\alpha,1}}(\overline{B_r(x)}) \geq |B_r(x)| = \omega_d r^d,$$

where $\omega_d$ is the volume of the unit ball. By monotonicity of the capacity, one deduces that

$$\text{cap}_{W^{\alpha,1}}(\{\mathcal{M}\varphi > t\}) \geq \text{cap}_{W^{\alpha,1}}(\overline{B_{\epsilon/t^{1/d}}(0)}) \geq \frac{\epsilon'}{t},$$

which implies (4.2) by integration with respect to $t$.

Although Proposition 4.1 fails for $\mathcal{M}u$, one does have a weak form of the capacitary inequality:

**Proposition 4.2.** *Let $\alpha > 0$. For every $u \in W^{\alpha,1}(\mathbb{R}^d)$ and $t > 0$, we have*

$$\text{cap}_{W^{\alpha,1}}(\{\mathcal{M}u > t\}) \leq \frac{C'}{t} \|u\|_{W^{\alpha,1}(\mathbb{R}^d)},$$

*for some constant $C' > 0$ depending on $\alpha$ and $d$.*

*Proof.* Given $t > 0$, observe that

(4.4) $$\{\mathcal{M}u > t\} \subset \{\mathcal{M}_1 u > t\} \cup A_t,$$

where

$$A_t := \{\mathcal{M}u > t \text{ and } \mathcal{M}_1 u \leq t\}.$$

By Proposition 4.1 and the Chebyshev inequality we have

$$\operatorname{cap}_{W^{\alpha,1}}\left(\{\mathcal{M}_1 u > t\}\right) \leq \frac{C}{t}\|u\|_{W^{\alpha,1}(\mathbb{R}^d)}.$$

By (4.4) and the subadditivity of the capacity, the proof is thus complete once we prove that

(4.5) $$\operatorname{cap}_{W^{\alpha,1}}(A_t) \leq \frac{C_1}{t}\|u\|_{L^1(\mathbb{R}^d)}.$$

The argument is based on the usual weak $L^1$ inequality for the maximal function. Indeed, using Wiener's covering lemma one finds sequences $(x_n)_{n\in\mathbb{N}}$ in $A_t$ and $(r_n)_{n\in\mathbb{N}}$ in the interval $(1,\infty)$ such that the balls $B_{r_n}(x_n)$ are disjoint, $A_t \subset \bigcup_{n\in\mathbb{N}} B_{5r_n}(x_n)$ and

$$t < \fint_{B_{r_n}(x_n)} |u|.$$

By (2.5) or (2.6), depending on $\alpha$, and the fact that $r_n > 1$,

$$\operatorname{cap}_{W^{\alpha,1}}(B_{5r_n}(x_n)) \leq C_2\,(5r_n)^d.$$

By countable subadditivity of the capacity and additivity of the integral, we thus have

$$\operatorname{cap}_{W^{\alpha,1}}(A_t) \leq \sum_{n=0}^{\infty} \operatorname{cap}_{W^{\alpha,1}}(B_{5r_n}(x_n)) \leq C_3 \sum_{n=0}^{\infty} r_n^d \leq \frac{C_4}{t} \int_{\bigcup_{n\in\mathbb{N}} B_{r_n}(x_n)} |u|,$$

which implies (4.5) and completes the proof. $\square$

## 5. Decomposition with $L^\infty$ bounds

We now show how one can obtain a decomposition in $W^{\alpha,1}(\mathbb{R}^d)$ that inherits $L^\infty$ bounds:

**Proposition 5.1.** *If, in addition to the assumptions of Theorems 1.1 or 3.1, we have $u \in L^\infty(\mathbb{R}^d)$, then the functions $u^\oplus$ and $u^\ominus$ can be chosen with the additional property that they also belong to $L^\infty(\mathbb{R}^d)$ and*

$$\|u^\oplus\|_{L^\infty(\mathbb{R}^d)} + \|u^\ominus\|_{L^\infty(\mathbb{R}^d)} \leq C'\|u\|_{L^\infty(\mathbb{R}^d)},$$

*for some constant $C' > 0$ depending on $\alpha$ and $d$.*

We rely on a clever construction of Harvey and Polking's (see Lemma 1 in [8]) that yields almost minimizers of the $W^{\alpha,1}$ capacity, with uniform bounds. Let us first illustrate this tool to get an improvement of Proposition 3.2:

**Proposition 5.2.** *For every $\varphi \in C_c^\infty(\mathbb{R}^N)$, one can find $\psi \in C_c^\infty(\mathbb{R}^N)$ satisfying, in addition to the conclusion of Proposition 3.2,*

$$\|\psi\|_{L^\infty(\mathbb{R}^d)} \leq 8\|\varphi\|_{L^\infty(\mathbb{R}^d)}.$$

To prove Proposition 5.2 we need the following

**Lemma 5.3.** *Let $\alpha > 0$. For every compact subset $K \subset \mathbb{R}^d$ and every $\epsilon > 0$, there exists $\varphi \in C_c^\infty(\mathbb{R}^d)$ such that $0 \leq \varphi \leq 1$ in $\mathbb{R}^d$, $\varphi = 1$ in a neighborhood of $K$, and*
$$\|\varphi\|_{W^{\alpha,1}(\mathbb{R}^d)} \leq C'' \operatorname{cap}_{W^{\alpha,1}}(K) + \epsilon,$$
*for some constant $C'' > 0$ depending on $\alpha$ and $d$.*

*Proof of Lemma 5.3.* We first assume that $\alpha < d$. Take a finite family $(Q_j)_{j\in\{0,\ldots,M\}}$ of dyadic closed cubes with disjoint interiors such that $K \subset \operatorname{int}\left(\bigcup_{j=0}^{M} Q_j\right)$, where $K \subset \mathbb{R}^d$ is a non-empty compact subset. By Lemma 1 in [8], there exist functions $(\varphi_j)_{j\in\{0,\ldots,M\}}$ with $\varphi_j \in C_c^\infty(\frac{3}{2}Q_j)$ such that

$$(5.1) \qquad 0 \leq \varphi_j \leq 1 \text{ in } \mathbb{R}^d, \qquad \sum_{j=0}^{M} \varphi_j = 1 \text{ in } \bigcup_{j=0}^{M} Q_j$$

and, for any $i \in \mathbb{N}$, the pointwise estimate holds

$$(5.2) \qquad |D^i \varphi_j| \leq \frac{C_1}{r_j^i} \quad \text{in } \mathbb{R}^d,$$

where $r_j > 0$ is the side length of $Q_j$ and the constant $C_1 > 0$ depends on $i$ and $d$, but not on the number $M+1$ of cubes. Since $\varphi_j$ is supported in $\frac{3}{2}Q_j$, one has in particular

$$(5.3) \qquad \|D^i \varphi_j\|_{L^1(\mathbb{R}^d)} \leq C_2\, r_j^{d-i}.$$

For every $0 < \theta < 1$, we then have by interpolation

$$(5.4) \qquad [D^i \varphi_j]_{W^{\theta,1}(\mathbb{R}^d)} \leq C_3 \|D^i \varphi_j\|_{L^1(\mathbb{R}^d)}^{1-\theta} \|D^{i+1}\varphi_j\|_{L^1(\mathbb{R}^d)}^{\theta} \leq C_4\, r_j^{d-i-\theta}.$$

To conclude the proof of the lemma, let $\delta > 0$ and assume that each side length $r_j$ is such that $r_j \leq \delta$. Hence, whether $\alpha$ is integer or not, it follows from (5.3) and (5.4) that $\varphi := \sum_{j=0}^{M} \varphi_j$ verifies

$$(5.5) \qquad \|\varphi\|_{W^{\alpha,1}(\mathbb{R}^d)} \leq C_5 \sum_{j=0}^{M} r_j^{d-\alpha}.$$

By the definition of the dyadic Hausdorff outer measure $\widehat{\mathcal{H}}_\delta^{d-\alpha}$ (see [14]), given $\epsilon > 0$ we take the cubes $(Q_j)_{j\in\{0,\ldots,M\}}$ so as to have

$$(5.6) \qquad \sum_{j=0}^{M} r_j^{d-\alpha} \leq \widehat{\mathcal{H}}_\delta^{d-\alpha}(K) + \epsilon.$$

Hence, by estimates (5.5) and (5.6) and the equivalence between $\widehat{\mathcal{H}}_\delta^{d-\alpha}$ and $\mathcal{H}_\delta^{d-\alpha}$,

$$\|\varphi\|_{W^{\alpha,1}(\mathbb{R}^d)} \leq C_6(\widehat{\mathcal{H}}_\delta^{d-\alpha}(K) + \epsilon) \leq C_7 \operatorname{cap}_{W^{\alpha,1}}(K) + C_6 \epsilon,$$

which proves the lemma for $\alpha < d$. When $\alpha \geq d$, one has from (5.5) that

$$\|\varphi\|_{W^{\alpha,1}(\mathbb{R}^d)} \leq C_8 \delta^{d-\alpha}(M+1).$$

Then, by minimization of the right-hand side with respect to $M$,

$$\|\varphi\|_{W^{\alpha,1}(\mathbb{R}^d)} \leq C_9 \, \mathcal{H}_\delta^0(K).$$

The conclusion thus follows in this case using Proposition 2.4. $\square$

*Proof of Proposition 5.2.* We rely on the construction from the proof of Proposition 3.2. By the previous lemma, we can pick each $\psi_j$ with the additional property that $0 \leq \psi_j \leq 2$ in $\mathbb{R}^d$. Therefore, taking $J$ as the smallest integer such that $|\varphi| \leq 2^J$ in $\mathbb{R}^d$, we get

$$0 \leq \psi = 2^m \zeta + \sum_{j=m}^{J-1} 2^{j+1} \psi_j \leq \sum_{j=-\infty}^{J-1} 2^{j+2} = 2^{J+2} \leq 8\|\varphi\|_{L^\infty(\mathbb{R}^d)},$$

which gives the $L^\infty$ bound for $\psi$. $\square$

In the proof of Proposition 5.1 we need the following counterpart of Lemma 5.3 on open sets:

**Lemma 5.4.** *Let $0 < \alpha < d$. For every open subset $\omega \subset \mathbb{R}^d$, there exists $v \in W^{\alpha,1}(\mathbb{R}^d)$ such that $0 \leq v \leq 1$ in $\mathbb{R}^d$, $v = 1$ in $\omega$, and*

$$\|v\|_{W^{\alpha,1}(\mathbb{R}^d)} \leq C'' \operatorname{cap}_{W^{\alpha,1}}(\omega).$$

*Proof of Lemma 5.4.* Let $0 < \delta < \infty$ and let $(Q_j)_{j \in \mathbb{N}}$ be a sequence of closed dyadic cubes with disjoint interiors and side lengths $r_j \leq \delta$ such that $\omega \subset \operatorname{int}\left(\bigcup_{j \in \mathbb{N}} Q_j\right)$ and

(5.7) $$\sum_{j=0}^{\infty} r_j^{d-\alpha} \leq 2\, \widehat{\mathcal{H}}_\delta^{d-\alpha}(\omega).$$

We rely on Harvey and Polking's lemma to construct a sequence of functions $(\varphi_j)_{j \in \mathbb{N}}$ with $\varphi_j \in C_c^\infty(\frac{3}{2}Q_j)$ that verifies properties (5.1)–(5.4) for every $j, M \in \mathbb{N}$. Some care is needed here because their lemma involves finitely many cubes, even though the constants in the estimates do not depend on the number of cubes. To deal with countably many cubes one may proceed as follows. Since $\sum_{j=0}^{\infty} r_j^{d-\alpha} < \infty$, the sequence of positive numbers $(r_j)_{j \in \mathbb{N}}$ converges to zero. We can thus relabel the cubes if necessary so that $(r_j)_{j \in \mathbb{N}}$ is non-increasing. Once the side lengths $r_j$ are arranged in this way and we have chosen the functions $\varphi_0, \ldots, \varphi_j$, the next function $\varphi_{j+1}$ is chosen in their proof without modification of the previous ones and for every $j \in \mathbb{N}$ we have

$$\|\varphi_j\|_{W^{\alpha,1}(\mathbb{R}^d)} \leq C_1 r_j^{d-\alpha}.$$

Hence, the series $v := \sum_{j=0}^{\infty} \varphi_j$ converges normally in $W^{\alpha,1}(\mathbb{R}^d)$ and in particular $v \in W^{\alpha,1}(\mathbb{R}^d)$. The conclusion follows from (5.7) and the equivalence between $\widehat{\mathcal{H}}_\delta^{d-\alpha}(\omega)$ and $\mathrm{cap}_{W^{\alpha,1}}(\omega)$. $\square$

In the regime $\alpha \notin \mathbb{N}$, Lemma 5.4 is a straightforward consequence of Lemma 5.3. Indeed, one takes a non-decreasing sequence of compact subsets $(K_j)_{j \in \mathbb{N}}$ whose capacities converge to $\mathrm{cap}_{W^{\alpha,1}}(\omega)$. For each $j \in \mathbb{N}$, by Lemma 5.3 there exists $\varphi_j \in C_c^\infty(\mathbb{R}^d)$ with $\varphi_j = 1$ on $K_j$ and

$$\|\varphi_j\|_{W^{\alpha,1}(\mathbb{R}^d)} \le C'' \mathrm{cap}_{W^{\alpha,1}}(K_j) + \frac{1}{2^j}.$$

For $\alpha = k + \theta$ with $k \in \mathbb{N}$ and $0 < \theta < 1$, by compactness one can extract a sequence $(\varphi_{j_i})_{i \in \mathbb{N}}$ that converges to some function $v$ in $W^{k,1}(\mathbb{R}^d)$. One has $v = 1$ in $\omega$, while Fatou's lemma implies that $D^k v \in W^{\theta,1}(\mathbb{R}^d)$. By lower semicontinuity of the norm, $v$ satisfies the required estimate. The adaptation of this argument to the case $\alpha = k \in \mathbb{N}_*$ is trickier since it yields a function $v$ whose distribution $D^k v$ could be no better than a finite measure in $\mathbb{R}^d$.

*Proof of Proposition 5.1.* Let $u \in (W^{\alpha,1} \cap L^\infty)(\mathbb{R}^d)$, which we assume in the first part of the proof to have compact support in some cube $Q \subset \mathbb{R}^d$. Given $\eta > 0$ to be explicitly chosen later on, take $\varphi \in C_c^\infty(\mathbb{R}^d)$ such that

(5.8) $\quad \|\varphi - u\|_{W^{\alpha,1}(\mathbb{R}^d)} \le \eta \quad \text{and} \quad \|\varphi\|_{L^\infty(\mathbb{R}^d)} \le \|u\|_{L^\infty(\mathbb{R}^d)}.$

For instance, $\varphi$ can be a convolution of $u$ with a smooth mollifier.

We begin by analyzing the case where $\alpha < d$. Given $\epsilon > 0$, we have

(5.9) $\quad\quad\quad |u - \varphi| \le \epsilon \quad \text{almost everywhere in } \mathbb{R}^d \setminus \omega,$

where $\omega$ is the open set $\{\mathcal{M}_1(\varphi - u) > \epsilon\}$. Let $v \in W^{\alpha,1}(\mathbb{R}^d)$ be as in Lemma 5.4. We then have

$$|u - \varphi| \le \epsilon + 2\|u\|_{L^\infty(\Omega)} v \quad \text{almost everywhere in } \mathbb{R}^d.$$

We now let $\psi \in C_c^\infty(\mathbb{R}^d)$ with $|\varphi| \le \psi$ in $\mathbb{R}^d$ given by Proposition 5.2. Then,

$$|u| \le |\varphi| + |u - \varphi| \le \psi + \epsilon + 2\|u\|_{L^\infty(\Omega)} v.$$

In addition,

$$\|\psi\|_{W^{\alpha,1}(\mathbb{R}^d)} \le C_1 \|\varphi\|_{W^{\alpha,1}(\mathbb{R}^d)} \le C_1(\|u\|_{W^{\alpha,1}(\mathbb{R}^d)} + \eta)$$

and, by choice of $v$ and the maximal strong capacitary inequality (Proposition 4.1),

$$\|v\|_{W^{\alpha,1}(\mathbb{R}^d)} \le C'' \mathrm{cap}_{W^{\alpha,1}}(\omega) \le \frac{C_2}{\epsilon} \|\varphi - u\|_{W^{\alpha,1}(\mathbb{R}^d)} \le \frac{C_2 \eta}{\epsilon}.$$

Take some fixed non-negative function $\zeta \in C_c^\infty(2Q)$ such that $0 \le \zeta \le 1$ in $\mathbb{R}^d$ and $\zeta = 1$ on $Q$ with uniform bounds on the derivatives depending on the side length of $Q$ and on the dimension $d$. Since $u$ is supported on $Q$,

$$|u| \le (\psi + \epsilon + 2\|u\|_{L^\infty(\Omega)} v)\zeta =: w. \tag{5.10}$$

This function $w$ belongs to $(W^{\alpha,1} \cap L^\infty)(\mathbb{R}^d)$, has compact support in $2Q$, and satisfies

$$\|w\|_{W^{\alpha,1}(\mathbb{R}^d)} \le C_3\big(\|\psi\|_{W^{\alpha,1}(\mathbb{R}^d)} + \epsilon + \|u\|_{L^\infty(\Omega)}\|v\|_{W^{\alpha,1}(\mathbb{R}^d)}\big)$$
$$\le C_4\Big(\|u\|_{W^{\alpha,1}(\mathbb{R}^d)} + \eta + \epsilon + \|u\|_{L^\infty(\Omega)}\frac{\eta}{\epsilon}\Big)$$

and

$$\|w\|_{L^\infty(\mathbb{R}^d)} \le 10\|u\|_{L^\infty(\Omega)} + \epsilon.$$

We now first choose $\epsilon > 0$ with

$$\epsilon \le \min\big\{\|u\|_{W^{\alpha,1}(\mathbb{R}^d)}, \|u\|_{L^\infty(\Omega)}\big\},$$

and then $\eta > 0$ so that

$$\eta + \|u\|_{L^\infty(\Omega)}\frac{\eta}{\epsilon} \le \|u\|_{W^{\alpha,1}(\mathbb{R}^d)}.$$

This concludes the proof when $u$ is compactly supported in a cube $Q$ and $\alpha < d$.

When $\alpha \ge d$, the function $u$ is continuous and, for any given $\epsilon > 0$, we can chose $\eta > 0$ sufficiently small so that (5.8) implies $\|\varphi - u\|_{L^\infty(\mathbb{R}^d)} \le \epsilon$. Thus, (5.9) holds in this case with $\omega = \emptyset$. The previous computation with $v = 0$ thus gives $w \in (W^{\alpha,1} \cap L^\infty)(\mathbb{R}^d)$ supported in $2Q$ such that (5.10) is satisfied and

$$\|w\|_{W^{\alpha,1}(\mathbb{R}^d)} \le C_4\big(\|u\|_{W^{\alpha,1}(\mathbb{R}^d)} + \eta + \epsilon\big) \quad \text{and} \quad \|w\|_{L^\infty(\mathbb{R}^d)} \le 8\|u\|_{L^\infty(\Omega)} + \epsilon.$$

The conclusion for $\alpha \ge d$ then follows from a suitable choice of $\epsilon$ and then $\eta$. The proof is thus complete for any $\alpha > 0$ when $u$ is supported in a cube.

For an arbitrary function $u \in (W^{\alpha,1} \cap L^\infty)(\mathbb{R}^d)$, we now decompose $\mathbb{R}^d$ as a countable union of closed cubes $(Q_j)_{j \in \mathbb{N}}$ with disjoint interiors and vertices given by all integer components. Take a sequence of functions $(\zeta_j)_{j \in \mathbb{N}}$ such that $\zeta_j$ is supported in $\frac{3}{2}Q_j$ and $\sum_{j=0}^\infty \zeta_j = 1$ in $\mathbb{R}^d$ obtained by translation of a single function $\zeta$. We then write $u = \sum_{j=0}^\infty u\zeta_j$. By Lemma 2.2,

$$\sum_{j=0}^\infty \|u\zeta_j\|_{W^{\alpha,1}(\mathbb{R}^d)} \le C_5\|u\|_{W^{\alpha,1}(\mathbb{R}^d)}.$$

We then apply the first part of the proof to each function $u\zeta_j$ to obtain a function $w_j \in (W^{\alpha,1} \cap L^\infty)(\mathbb{R}^d)$ supported in $2Q_j$. Since the cubes $2Q_j$ overlap a finite number of times, the function $u^\oplus := \sum_{j=0}^\infty w_j$ also belongs to $(W^{\alpha,1} \cap L^\infty)(\mathbb{R}^d)$. To complete the proof, it suffices to take $u^\ominus := u^\oplus - u$. $\square$

ACKNOWLEDGEMENTS

This work was initiated while the first author (ACP) was visiting the Department of Applied Mathematics of the National Chiao Tung University with support from the National Center for Theoretical Sciences, Wallonie-Bruxelles International and the Fonds de la Recherche scientifique–FNRS under research grant J.0020.18. He warmly thanks the Department for the invitation and hospitality. The second author (DS) is supported by the Taiwan Ministry of Science and Technology under research grants 105-2115-M-009-004-MY2 and 107-2918-I-009-003.

REFERENCES

[1] D. R. Adams, *On the existence of capacitary strong type estimates in $\mathbb{R}^n$*, Ark. Mat. **14** (1976), 125–140. ↑3

[2] ______ , *A note on Choquet integrals with respect to Hausdorff capacity*, Function spaces and applications (Lund, 1986), Lecture Notes in Math., vol. 1302, Springer, Berlin, 1988, pp. 115–124. ↑4, 10

[3] D. R. Adams and L. I. Hedberg, *Function spaces and potential theory*, Grundlehren der Mathematischen Wissenschaften, vol. 314, Springer-Verlag, Berlin, 1996. ↑2

[4] A. Carlsson and V. G. Maz'ya, *On approximation in weighted Sobolev spaces and self-adjointness*, Math. Scand. **74** (1994), 111–124. ↑6

[5] B. E. J. Dahlberg, *A note on Sobolev spaces*, Harmonic analysis in Euclidean spaces (Proc. Sympos. Pure Math., Williams Coll., Williamstown, MA, 1978), Part 1, Proc. Sympos. Pure Math., XXXV, Amer. Math. Soc., Providence, RI, 1979, pp. 183–185. ↑2

[6] H. Federer and W. P. Ziemer, *The Lebesgue set of a function whose distribution derivatives are p-th power summable*, Indiana Univ. Math. J. **22** (1972/73). ↑2

[7] W. Gustin, *Boxing inequalities*, J. Math. Mech. **9** (1960), 229–239. ↑2

[8] R. Harvey and J. Polking, *Removable singularities of solutions of linear partial differential equations*, Acta Math. **125** (1970), 39–56. ↑3, 12, 13

[9] V. G. Maz'ya, *Classes of sets and measures that are connected with imbedding theorems*, Imbedding theorems and their applications (Proc. Sympos., Baku, 1966), Izdat. "Nauka", Moscow, 1970, pp. 142–159 (Russian). ↑3

[10] ______ , *Strong capacity-estimates for "fractional" norms*, Zap. Naučn. Sem. Leningrad. Otdel. Mat. Inst. Steklov. (LOMI) **70** (1977), 161–168, 292 (Russian); English transl., J. Soviet Math. **23** (1983), 1997–2003. ↑3

[11] V. G. Maz'ya and T. O. Shaposhnikova, *Theory of Sobolev multipliers*, Grundlehren der Mathematischen Wissenschaften, vol. 337, Springer-Verlag, Berlin, 2009. ↑3

[12] J. Orobitg and J. Verdera, *Choquet integrals, Hausdorff content and the Hardy-Littlewood maximal operator*, Bull. London Math. Soc. **30** (1998), 145–150. ↑10

[13] A. C. Ponce and D. Spector, *A boxing inequality for the fractional perimeter*. To appear in Ann. Scuola Norm. Sup. Pisa Cl. Sci. ↑3, 4, 6, 8, 10

[14] D. Yang and W. Yuan, *A note on dyadic Hausdorff capacities*, Bull. Sci. Math. **132** (2008), 500–509. ↑13

[15] W. P. Ziemer, *Weakly differentiable functions*, Graduate Texts in Mathematics, vol. 120, Springer-Verlag, New York, 1989. ↑2

AUGUSTO C. PONCE
UNIVERSITÉ CATHOLIQUE DE LOUVAIN
INSTITUT DE RECHERCHE EN MATHÉMATIQUE ET PHYSIQUE
CHEMIN DU CYCLOTRON 2, L7.01.02
1348 LOUVAIN-LA-NEUVE
BELGIUM
*E-mail address*: Augusto.Ponce@uclouvain.be

DANIEL SPECTOR
NATIONAL CHIAO TUNG UNIVERSITY
DEPARTMENT OF APPLIED MATHEMATICS
HSINCHU, TAIWAN

NATIONAL CENTER FOR THEORETICAL SCIENCES
NATIONAL TAIWAN UNIVERSITY
NO. 1 SEC. 4 ROOSEVELT RD.
TAIPEI, 106, TAIWAN
*E-mail address*: dspector@math.nctu.edu.tw